\documentclass[12pt]{article}
\usepackage{latexsym,amssymb,amsmath,enumerate,float,geometry,cite,xcolor,graphicx}
\newtheorem{theorem}{Theorem}

\newtheorem{lemma}[theorem]{Lemma}
\newtheorem{conjecture}[theorem]{Conjecture}

\begin{document}


\title{\Large On the K\H{o}nig-Egerv\'ary Theorem for $k$-Paths}
\author{St\'{e}phane Bessy$^1$, Pascal Ochem$^1$, and Dieter Rautenbach$^2$}
\date{}
\maketitle
\vspace{-10mm}
\begin{center}
{\small
$^1$ 
Laboratoire d'Informatique, de Robotique et de Micro\'{e}lectronique de Montpellier,\\
Montpellier, France, \texttt{stephane.bessy@lirmm.fr,pascal.ochem@lirmm.fr}\\[3mm]
$^2$ Institute of Optimization and Operations Research, Ulm University,\\
Ulm, Germany, \texttt{dieter.rautenbach@uni-ulm.de}}
\end{center}

\begin{abstract}
The famous K\H{o}nig-Egerv\'ary theorem is equivalent to the statement 
that the matching number equals the vertex cover number 
for every induced subgraph of some graph
if and only if that graph is bipartite.
Inspired by this result, we consider the set ${\cal G}_k$ of all graphs such that, 
for every induced subgraph,
the maximum number of disjoint paths of order $k$
equals the minimum order of a set of vertices intersecting all paths of order $k$.
For $k\in \{ 3,4\}$, we give complete structural descriptions of the graphs in ${\cal G}_k$.
Furthermore, for odd $k$, we give a complete structural description of the graphs in ${\cal G}_k$
that contain no cycle of order less than $k$.
For these graph classes,
our results yield efficient recognition algorithms as well as 
efficient algorithms that determine maximum sets of disjoint paths of order $k$
and minimum sets of vertices intersecting all paths of order $k$.

\end{abstract}

{\small
\begin{tabular}{lp{12.5cm}}
\textbf{Keywords:} & K\H{o}nig-Egerv\'ary theorem; matching; vertex cover; $k$-path vertex cover; bipartite graph\\
\textbf{MSC2010:} & 05C69
\end{tabular}
}

\pagebreak

\section{Introduction}

The famous K\H{o}nig-Egerv\'ary theorem \cite{k,e} states that 
the matching number $\nu(G)$ of a bipartite graph $G$ equals 
its vertex cover number $\tau(G)$. 
Since a graph is bipartite if and only if it contains no odd cycle $C_{2k+1}$ as an induced subgraph, 
and $\nu(C_{2k+1})=k<k+1=\tau(C_{2k+1})$,
the K\H{o}nig-Egerv\'ary theorem is equivalent to the statement 
that $\nu(H)=\tau(H)$ for every induced subgraph $H$ of some graph $G$
if and only if $G$ is bipartite.
Considering a matching as a packing of paths of order $2$,
and a vertex cover as a set of vertices intersecting every path of order $2$, 
it is natural to ask for generalizations of the K\H{o}nig-Egerv\'ary theorem for longer paths,
and to consider the corresponding graph classes generalizing the bipartite graphs.

In the present paper we study such generalizations.

We consider finite, simple, and undirected graphs
as well as finite and undirected multigraphs 
that may contain loops and parallel edges.
Let $k$ be a positive integer, and let $G$ be a graph.
A {\it $k$-path} and a {\it $k$-cycle} in $G$ is a not necessarily induced path and cycle of order $k$ in $G$, respectively.
A set of disjoint $k$-paths in $G$ 
is a {\it $k$-matching} in $G$,
and a set of vertices of $G$ intersecting every $k$-path in $G$ 
is a {\it $k$-vertex cover} in $G$.
The {\it $k$-matching number} $\nu_k(G)$ of $G$ is the maximum cardinality of a $k$-matching in $G$,
and the {\it $k$-vertex cover number} $\tau_k(G)$ of $G$
is the minimum cardinality of a $k$-vertex cover in $G$.

Clearly, $$\nu_k(G)\leq\tau_k(G).$$
Let ${\cal G}_k$ be the set of all graphs $G$
such that $\nu_k(H)=\tau_k(H)$ for every induced subgraph $H$ of $G$. 
As noted above, the K\H{o}nig-Egerv\'ary theorem is equivalent to the statement 
that ${\cal G}_2$ is the set of all bipartite graphs.
Since $\nu_1(G)=\tau_1(G)=n(G)$ for every graph $G$ of order $n(G)$,
the set ${\cal G}_1$ contains all graphs.

For $k\in \{ 3,4\}$, we give complete structural descriptions of the graphs in ${\cal G}_k$.
Furthermore, for odd $k$, we give a complete structural description of the graphs in ${\cal G}_k$
that contain no cycle of order less than $k$.

Among the two parameters $\nu_k(G)$ and $\tau_k(G)$, 
only the latter seems to have received considerable attention in the literature \cite{bjkst,bkks,lzh}.
Note that a set $X$ of vertices of a graph $G$ 
is a $3$-vertex cover if and only if its complement $V(G)\setminus X$ 
is a so-called {\it dissociation set} \cite{y,bcl}, that is, 
a set of vertices inducing a subgraph of maximum degree at most $1$.
Probably motivated by this connection,
the $3$-vertex cover number has been studied in detail \cite{kks,k2,tz,xk,twyc}.
For every $k$ at least $3$,
the hardness of the $k$-vertex cover number has been shown in \cite{bkks}. 
It follows from known results (cf. [GT12] in \cite{gj})
that, for every integer $k$ at least $3$, 
it is NP-complete to decide for a given graph $G$ 
whose order $n(G)$ is a multiple of $k$, 
whether $\nu_k(G)=\frac{n(G)}{k}$,
that is, whether $G$ has a {\it perfect} $k$-matching.

For a positive integer $k$ and a graph $G$, 
let ${\cal P}_k$ be the set of all $k$-paths of $G$.
The parameters $\nu_k(G)$ and $\tau_k(G)$ are the optimum values of the following integer linear programs.
$$
\nu_k(G)\left\{\hspace{1cm}\begin{array}{rrcll}
\max &	\sum\limits_{P\in {\cal P}_k}x_P &  & & \\[5mm]
s.t. & \sum\limits_{P\in {\cal P}_k:\,\, u\in V(P)}x_P & \leq & 1 & \forall u\in V(G)\\
&	x_P & \in & \{ 0,1\} & \forall P\in {\cal P}_k
\end{array}\right.$$

$$\tau_k(G)\left\{\,\,\,\,\,\,\,\,\,\,\,\,\,\,\,\,\hspace{1cm}\begin{array}{rrcll}
\min &	\sum\limits_{u\in V(G)}y_u &  & & \\[5mm]
s.t. & \sum\limits_{u\in V(P)}y_u & \geq & 1 & \forall P\in {\cal P}_k\\
&	y_u & \in & \{ 0,1\} & \forall u\in V(G)
\end{array}\right.
$$
Relaxing ``$\in \{ 0,1\}$'' in both programs to ``$\geq 0$'' yields a pair of dual linear programs,
whose optimal values we denote by $\nu_k^*(G)$ and $\tau_k^*(G)$, respectively.
Since $\nu_k(G)=\nu_k^*(G)=\tau_k^*(G)=\tau_k(G)$ for a given graph $G$ in ${\cal G}_k$, 
linear programming allows to determine $\nu_k(G)$ and $\tau_k(G)$ for $G$ in polynomial time.
Furthermore, 
since ${\cal G}_k$ is closed under taking induced subgraphs, 
iteratively considering the removal of individual vertices, 
one can use linear programming to determine in polynomial time 
an induced subgraph $G'$ of $G$ of minimum order 
with $\nu_k(G)=\nu_k(G')=\tau_k(G')=\tau_k(G)$.
Note that a maximum $k$-matching in $G'$ covers all vertices of $G'$, and is also a maximum $k$-matching in $G$,
and that a minimum $k$-vertex cover in $G'$ is also a minimum $k$-vertex cover in $G$.
Now, within $G'$, one can use linear programming to 
iteratively identify in polynomial time $k$-paths as well as vertices 
whose removal reduces the $k$-matching number as well as $k$-vertex cover
by exactly $1$, respectively.
Clearly, the identified $k$-paths form a maximum $k$-matching in $G$, 
and the identified vertices form a minimum $k$-vertex cover in $G$.

We discuss some generic examples of graphs in ${\cal G}_k$, namely, 
\begin{itemize}
\item forests, 
\item $k$-subdivisions of multigraphs, and, 
\item $k/2$-subdivisions of bipartite multigraphs for even $k$.
\end{itemize}
Trivially, every graph of order less than $k$ belongs to ${\cal G}_k$,
which implies that the local structure of the graphs in ${\cal G}_k$ is not simple.

The fact that all forests belong to ${\cal G}_k$ follows by a inductive argument using the following lemma.
In fact, the lemma yields a simple polynomial time reduction algorithm that determines
a maximum $k$-matching as well as a minimum $k$-vertex cover in a given forest. 
An efficient algorithm computing a minimum $k$-vertex cover in a given forest was presented in \cite{bkks}.

\begin{lemma}\label{lemmatree}
Let $k$ be a positive integer.
If the graph $G$ is the union of a tree $T$ and a graph $G'$ such that $T$ and $G'$ share exactly one vertex $x$, 
the tree $T$ contains a $k$-path, 
but the forest $T-x$ contains no $k$-path,
then 
$\nu_k(G)=\nu_k(G'-x)+1$ and $\tau_k(G)=\tau_k(G'-x)+1$.
\end{lemma}
{\it Proof:} Every $k$-path in $T$ contains $x$.
Hence, if ${\cal P}$ is a $k$-matching in $G$, 
then at most one path in ${\cal P}$ intersects $V(T)$.
Removing any such path yields a $k$-matching in $G'-x$, 
which implies $\nu_k(G)\leq \nu_k(G'-x)+1$.
Conversely, if ${\cal P}'$ is a $k$-matching in $G'-x$, 
then adding a $k$-path contained in $T$, 
yields a $k$-matching in $G$,
which implies $\nu_k(G)\geq \nu_k(G'-x)+1$.

If $X$ is a $k$-vertex cover in $G$, then $X$ intersects $V(T)$,
and $X\setminus V(T)$ is a $k$-vertex cover in $G'-x$,
which implies $\tau_k(G)\geq \tau_k(G'-x)+1$.
Conversely, adding $x$ to any $k$-vertex cover in $G'-x$
yields a $k$-vertex cover in $G$,
which implies $\tau_k(G)\leq \tau_k(G'-x)+1$. $\Box$

\medskip

\noindent The following lemma captures some natural cycle conditions for the graphs in ${\cal G}_k$.

For an integer $n$, let $[n]$ be the set of positive integers at most $n$.

\begin{lemma}\label{lemmacycle}
Let $k$ and $p$ be positive integers.
\begin{enumerate}[(i)]
\item Every cycle of order at least $k$ in every graph in ${\cal G}_k$ has order $0$ modulo $k$.
\item A set $X$ of vertices of the cycle $C_{pk}:u_1u_2\ldots u_{pk}u_1$ of order $pk$ 
is a minimum $k$-vertex cover in $C_{pk}$ if and only if 
$X=\{ u_{i+(j-1)k}:j\in [p]\}$ for some $i\in [k]$.
\item If $G$ is in ${\cal G}_3$, 
$C$ is a cycle in $G$,
and $u$ and $v$ are distinct vertices of $C$ that have neighbors outside of $V(C)$, then ${\rm dist}_C(u,v)\equiv 0$ {\rm mod} $3$.
\item If $G$ is in ${\cal G}_4$, 
$C$ is a cycle of length at least $4$ in $G$,
and $u$ and $v$ are distinct vertices of $C$ that have neighbors outside of $V(C)$, then ${\rm dist}_C(u,v)\equiv 0$ {\rm mod} $2$.
\end{enumerate}
\end{lemma}
{\it Proof:} If the graph $G$ arises by adding some edges to the cycle $C_n$ of order $n$, where $n$ is at least $k$, 
then $\nu_k(G)=\left\lfloor\frac{n}{k}\right\rfloor\leq\left\lceil\frac{n}{k}\right\rceil=\tau_k(C_n)\leq \tau_k(G)$,
which implies (i).
The value of $p=\tau_k(C_{pk})$ and the fact that every $k$-vertex cover in $C_{pk}$ 
has to contain at least one of any $k$ consecutive vertices of $C_{pk}$ implies (ii).

If $G$, $C$, $u$, and $v$ are as in (iii), 
$u'$ is a neighbor of $u$ outside of $V(C)$,
$v'$ is a neighbor of $v$ outside of $V(C)$,
and $G'$ is the subgraph of $G$ induced by $V(C)\cup \{ u',v'\}$,
then $\nu_3(G')=\left\lfloor\frac{n(C)+|\{ u',v'\}|}{3}\right\rfloor=\frac{n(C)}{3}$.
Since $G\in {\cal G}_3$, we obtain $\tau_3(G')=\frac{n(C)}{3}=\tau_3(C)$,
which implies that every minimum $3$-vertex cover in $G'$ is a minimum $3$-vertex cover in $C$,
and, hence, as described in (ii).
Since $u$ and $v$ must both belong to every minimum $3$-vertex cover in $G'$,
their distance on $C$ must be a multiple of $3$.

Now, if $G$, $C$, $u$, and $v$ are as in (iv),
and $u'$, $v'$, and $G'$ are as above,
then $\nu_4(G')=\left\lfloor\frac{n(C)+|\{ u',v'\}|}{4}\right\rfloor=\frac{n(C)}{4}$.
Again every minimum $4$-vertex cover in $G'$ is a minimum $4$-vertex cover in $C$,
and, hence, as described in (ii).
Since every minimum $4$-vertex cover in $G'$ contains either $u$ or both vertices at distance $2$ from $u$ within $C$,
and the same holds for $v$, the distance of $u$ and $v$ on $C$ must be even.
$\Box$

\medskip

\noindent Lemma \ref{lemmacycle} (i) and (iii) suggest that subdividing every edge of a multigraph $k-1$ times 
yields a natural candidate for a graph in ${\cal G}_k$.
For a positive integer $k$, 
let the {\it $k$-subdivision} $Sub_k(H)$ of a multigraph $H$ arise by subdividing every edge of $H$ exactly $k-1$ times, 
that is, 
\begin{itemize}
\item every edge between distinct vertices $u$ and $v$
is replaced by a $(k+1)$-path between $u$ and $v$ whose internal vertices have degree $2$, and 
\item every loop incident with some vertex $u$
is replaced by a $k$-cycle containing $u$ and $k-1$ further vertices of degree $2$.
\end{itemize}
Note that the $k$-subdivision of a forest is a forest.
Together with Lemma \ref{lemmatree}, the following lemma implies 
that $Sub_k(H)$ belongs to ${\cal G}_k$ for every multigraph $H$.

\begin{lemma}\label{lemmasub}
Let $k$ be a positive integer.
If the graph $G$ contains an induced subgraph $B$ such that 
\begin{itemize}
\item $B=Sub_k(H)$ for some connected multigraph $H$ that contains a cycle, and
\item every component $K$ of $G-V(H)$ that contains a vertex from $V(B)\setminus V(H)$ satisfies $\nu_k(K)=0$,
\end{itemize}
then $\nu_k(G)=\nu_k(G-V(H))+n(H)$, and $\tau_k(G)=\tau_k(G-V(H))+n(H)$.
\end{lemma}
{\it Proof:} Since $H$ is connected and contains a cycle, 
it contains an edge $e$ incident with some vertex $r$ such that $H-e$ contains a spanning tree $T$ of $H$.
Rooting $T$ in $r$, 
assigning $e$ to $r$, and 
assigning to every other vertex of $H$, 
the edge to its parent within $T$,
yields an injective function $f:V(H)\to E(H)$ such that $u$ is incident with $f(u)$ for every vertex $u$ of $H$.

Let ${\cal P}_f$ be $k$-matching of order $n(H)$ in $B$ that contains, for every vertex $u$ of $H$,
the $k$-path formed within $B$ by $u$ and the subdivided edge $f(u)$.
Recall that the components of $G-V(H)$ that contain a vertex from $V(B)\setminus V(H)$ 
contain no $k$-paths.
Therefore, adding ${\cal P}_f$ to any $k$-matching in $G-V(H)$ 
yields $\nu_k(G)\geq\nu_k(G-V(H))+n(H)$.
Conversely, if ${\cal P}$ is a $k$-matching in $G$,
then, since every $k$-path in $G$ that intersects $V(B)$ contains a vertex of $H$,
the set ${\cal P}$ contains at most $n(H)$ paths intersecting $V(B)$.
Removing all such paths from ${\cal P}$ yields a $k$-matching in $G-V(H)$,
which implies $\nu_k(G)\leq\nu_k(G-V(H))+n(H)$.

If $X$ is a $k$-vertex cover in $G-V(H)$, then $X\cup V(H)$ is a $k$-vertex cover in $G$,
which implies $\tau_k(G)\leq\tau_k(G-V(H))+n(H)$.
Now, let $X$ be a $k$-vertex cover in $G$.
Clearly, $X'=X\cap V(B)$ is a $k$-vertex cover in $B$.
If some vertex $u$ of $H$ does not belong to $X'$,
then $X'$ must intersect all subdivided edges of $H$ incident with $u$,
in particular, $X'$ contains a vertex from the subdivided edge $f(u)$.
Since $f$ is injective, 
this easily implies that $X'$ contains at least $n(H)$ vertices.
Since $X\setminus X'$ is a $k$-vertex cover in $G-V(H)$,
we obtain $\tau_k(G)\geq\tau_k(G-V(H))+n(H)$.
$\Box$

\medskip

\noindent For even values of $k$, Lemma \ref{lemmacycle} (i) and (iv) suggest yet another construction
based on subdivisions of bipartite multigraphs.
The following lemma captures the essence of this construction.

\begin{lemma}\label{lemmabipsub}
If $k$ is a positive even integer,
and $G=Sub_{k/2}(H)$ for some bipartite connected multigraph $H$ that contains a cycle,
then $\nu_k(G)=\tau_k(G)$.
\end{lemma}
{\it Proof:} 
In view of the K\H{o}nig-Egerv\'ary theorem, and, since $H$ is bipartite, 
it suffices to show that $\nu_k(G)\geq \nu(H)$ and $\tau_k(G)\leq \tau(H)$.

let $M$ be a matching in $H$.
Contracting the edges in $M$ yields a connected multigraph that contains a cycle,
and arguing similarly as in the proof of Lemma \ref{lemmasub},
we obtain the existence of an injective function $f:M\to E(H)\setminus M$ 
such that the edges $e$ and $f(e)$ are adjacent for every edge $e$ in $M$.
Now, for every edge $e$ in $M$,
the $(k/2+1)$-path corresponding to the subdivided edge $e$
and the $(k/2-1)$-path corresponding to the interior of the subdivided edge $f(e)$
form a $k$-path in $G$.
Since $M$ is a matching and $f$ is injective, all these $k$-paths are disjoint, 
which implies $\nu(H)\leq \nu_k(G)$.

If $X$ is a vertex cover in $H$, then every component of $G-X$ is a
$(k/2-1)$-subdivision of some star.
Hence, $G-X$ contains no $k$-path, which implies $\tau(H)\geq \tau_k(G)$.
$\Box$

\section{The graphs in ${\cal G}_3$ and ${\cal G}_4$}

In this section we characterize the graphs in ${\cal G}_k$ for $k\in \{ 3,4\}$
by describing their blocks and conditions imposed on their cutvertices.
As it turns out, the three generic examples of graphs in ${\cal G}_k$ discussed in the introduction
are the main building blocks of the considered graphs.

Recall that a cutvertex of a graph $G$ is a vertex $x$ of $G$ for which $G-x$ has more components than $G$,
and that a block of $G$ is a maximal connected subgraph $B$ of $G$ such that $B$ itself has no cutvertex.
An endblock of $G$ is a block of $G$ 
that contains at most one cutvertex of $G$.
A block is trivial if it is either $K_1$ or $K_2$.

Let ${\cal H}_3$ be the set of all graphs $G$ such that every non-trivial block $B$ of $G$ 
satisfies the following condition.
\begin{enumerate}[(i)]
\item $B=Sub_3(H)$ for some multigraph $H$, and every cutvertex of $G$ that belongs to $B$ is a vertex of $H$.
\end{enumerate}

\begin{theorem}\label{theoremg3}
${\cal G}_3={\cal H}_3$.
\end{theorem}
{\it Proof:} In order to show that ${\cal G}_3\subseteq {\cal H}_3$,
it suffices to show that $G\in {\cal H}_3$ for every connected graph $G$ in ${\cal G}_3$.
If $G$ is a tree, then all blocks of $G$ are trivial, and, hence, $G\in {\cal H}_3$. 
If $G$ is a cycle, then Lemma \ref{lemmacycle}(i) implies that $n(G)$ is a multiple of $3$,
and, hence, $G=Sub_3(C_{n(G)/3})\in {\cal H}_3$.
Now, we may assume that $G$ is neither a tree nor a cycle.
Let $B$ be a non-trivial block of $G$.
By Lemma \ref{lemmacycle}(i),
the order of every cycle in $B$ is a multiple of $3$.
Suppose that $B$ contains a path $P:u_0\ldots u_\ell$ 
such that $u_0$ and $u_\ell$ have degree at least $3$ in $G$,
and $u_1,\ldots,u_{\ell-1}$ have degree $2$ in $G$.
Since $B-u_1$ is connected, 
the path $P$ is contained in a cycle $C$ 
such that $u_0$ and $u_\ell$ both have neighbors outside of $V(C)$.
By Lemma \ref{lemmacycle}(iii), the length $\ell$ of $P$ is a multiple of $3$,
in particular, no two vertices of $B$ of degree at least $3$ in $G$ are adjacent.
Let $H$ be the multigraph that arises by replacing every path or cycle
$u_0u_1u_2u_3\ldots u_{3p-3}u_{3p-2}u_{3p-1}u_{3p}$ of length $3p$ 
such that $u_0$ and $u_{3p}$ have degree at least $3$ in $G$,
and $u_1,\ldots,u_{3p-1}$ have degree $2$ in $G$,
by the path or cycle $u_0u_3\ldots u_{3p-3}u_{3p}$ of length $p$. 
Clearly, $B=Sub_3(H)$, and every cutvertex of $G$ that belongs to $B$ is a vertex of $H$,
that is, $G\in {\cal H}_3$.
Altogether, we obtain ${\cal G}_3\subseteq {\cal H}_3$.

It follows easily from its definition that ${\cal H}_3$ is a hereditary class of graphs,
that is, it is closed under taking induced subgraphs.
Therefore, in order to show the reverse inclusion ${\cal H}_3\subseteq {\cal G}_3$,
it suffices to show that $\nu_3(G)=\tau_3(G)$ for every connected graph $G$ in ${\cal H}_3$,
which we do by induction on the order of $G$.
If $G$ is a tree, then Lemma \ref{lemmatree} implies $\nu_3(G)=\tau_3(G)$.
If $G$ is a cycle, then the order of $G$ is a multiple of $3$, and, hence, $\nu_3(G)=\tau_3(G)$.
Now, we may assume that $G$ is neither a tree nor a cycle.
Let $B$ be a non-trivial block of $G$.
Let $B=Sub_3(H)$ for some multigraph $H$ such that every cutvertex of $G$ that belongs to $B$ is a vertex of $H$.
By Lemma \ref{lemmasub} applied to $B$,
we obtain $\nu_3(G)=\nu_3(G-V(H))+n(H)$ and $\tau_3(G)=\tau_3(G-V(H))+n(H)$.
Since ${\cal H}_3$ is hereditary,
we obtain, by induction, $\nu_3(G-V(H))=\tau_3(G-V(H))$,
which implies $\nu_3(G)=\tau_3(G)$ and completes the proof. $\Box$

\medskip 

\noindent For some positive integer $p$, 
let the graph $T(p)$ arise by adding an edge between the two vertices in a partite set of order $2$ of the complete bipartite graph $K_{2,p}$.
Note that $T(1)$ is a triangle, and that $T(2)$ arises by removing one edge from $K_4$.

\begin{figure}[H]
\begin{center}
\unitlength 1mm 
\linethickness{0.4pt}
\ifx\plotpoint\undefined\newsavebox{\plotpoint}\fi 
\begin{picture}(38,44)(0,0)
\put(8,10){\circle*{2}}
\put(28,10){\circle*{2}}
\put(18,20){\circle*{2}}
\put(18,30){\makebox(0,0)[cc]{$\vdots$}}
\put(18,40){\circle*{2}}
\put(18,30){\oval(8,28)[]}
\put(8,10){\line(1,0){20}}
\put(28,10){\line(-1,1){10}}
\put(18,20){\line(-1,-1){10}}
\put(8,10){\line(1,3){10}}
\put(18,40){\line(1,-3){10}}
\qbezier(28,0)(15,23)(38,10)
\put(16,0){\circle*{2}}
\put(0,0){\circle*{2}}
\put(8,0){\makebox(0,0)[cc]{$\dots$}}
\put(16,0){\line(-4,5){8}}
\put(8,10){\line(-4,-5){8}}
\put(4,12){\makebox(0,0)[cc]{$x$}}
\end{picture}
\end{center}
\caption{$T(p)$ as a non-endblock of a graph in ${\cal G}_4$.}\label{fig1}
\end{figure}
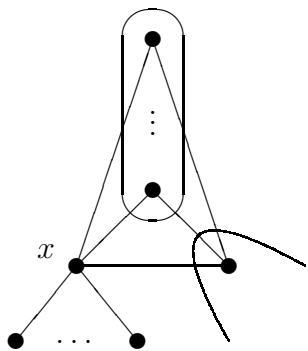

Let ${\cal H}_4$ be the set of all graphs $G$ such that every non-trivial block $B$ of $G$ 
satisfies the following condition.
\begin{enumerate}[(i)]
\item Either $B=Sub_2(H)$ for some bipartite multigraph $H$, and every cutvertex of $G$ that belongs to $B$ is a vertex of $H$,
\item or $B=K_4$ is an endblock,
\item or $B=T(2)$ is an endblock, and, if $B$ contains a cutvertex $x$ of $G$, then $x$ has degree $2$ in $B$,
\item or $B=T(p)$ for some positive integer $p$, 
at most two cutvertices of $G$ belong to $B$,
every cutvertex of $G$ that belongs to $B$ has degree $p+1$ in $B$,
and, if $B$ contains two cutvertices of $G$,
then there is one cutvertex $x$ of $G$ in $B$ such that every vertex in $N_G(x)\setminus V(B)$ has degree $1$ in $G$.
\end{enumerate}
See Figure \ref{fig1} for an illustration of (iv).

\begin{theorem}\label{theoremg4}
${\cal G}_4={\cal H}_4$.
\end{theorem}
{\it Proof:} As before, in order to show that ${\cal G}_4\subseteq {\cal H}_4$,
we show that $G\in {\cal H}_4$ for every connected graph $G\in {\cal G}_4$.
If $G$ is a tree, then clearly $G\in {\cal H}_4$.
If $G$ is a cycle, then Lemma \ref{lemmacycle}(i) implies that $n(G)$ is either $3$ or a multiple of $4$,
and, hence, $G\in {\cal H}_4$.
Now, we may assume that $G$ is neither a tree nor a cycle.
Let $B$ be a non-trivial block of $G$.

\begin{figure}[H]
\begin{center}
$\mbox{}$\hfill
\unitlength 0.8mm 
\linethickness{0.4pt}
\ifx\plotpoint\undefined\newsavebox{\plotpoint}\fi 
\begin{picture}(21,36)(0,0)
\put(0,15){\circle*{2}}
\put(20,15){\circle*{2}}
\put(10,25){\circle*{2}}
\put(0,15){\line(1,0){20}}
\put(20,15){\line(-1,1){10}}
\put(10,25){\line(-1,-1){10}}
\put(20,5){\circle*{2}}
\put(0,5){\circle*{2}}
\put(10,35){\circle*{2}}
\put(10,35){\line(0,-1){10}}
\put(0,15){\line(0,-1){10}}
\put(20,15){\line(0,-1){10}}
\put(10,0){\makebox(0,0)[cc]{$G_1$}}
\end{picture}
\hfill
\linethickness{0.4pt}
\ifx\plotpoint\undefined\newsavebox{\plotpoint}\fi 
\begin{picture}(21,36)(0,0)
\put(0,25){\circle*{2}}
\put(20,25){\circle*{2}}
\put(10,35){\circle*{2}}
\put(0,25){\line(1,0){20}}
\put(20,25){\line(-1,1){10}}
\put(10,35){\line(-1,-1){10}}
\put(20,15){\circle*{2}}
\put(0,15){\circle*{2}}
\put(0,25){\line(0,-1){10}}
\put(20,25){\line(0,-1){10}}
\put(10,0){\makebox(0,0)[cc]{$G_2$}}
\put(20,5){\circle*{2}}
\put(0,5){\circle*{2}}
\put(20,15){\line(0,-1){10}}
\put(0,15){\line(0,-1){10}}
\end{picture}
\hfill
\linethickness{0.4pt}
\ifx\plotpoint\undefined\newsavebox{\plotpoint}\fi 
\begin{picture}(21,46)(0,0)
\put(0,25){\circle*{2}}
\put(20,25){\circle*{2}}
\put(10,35){\circle*{2}}
\put(0,25){\line(1,0){20}}
\put(20,25){\line(-1,1){10}}
\put(10,35){\line(-1,-1){10}}
\put(10,0){\makebox(0,0)[cc]{$G_3$}}
\put(10,15){\circle*{2}}
\put(10,5){\circle*{2}}
\put(10,45){\circle*{2}}
\put(10,45){\line(0,-1){10}}
\put(20,25){\line(-1,-1){10}}
\put(10,15){\line(-1,1){10}}
\put(10,15){\line(0,-1){9}}
\end{picture}
\hfill$\mbox{}$
\end{center}
\caption{Three forbidden subgraphs for the graphs in ${\cal G}_4$.}\label{fig2}
\end{figure}
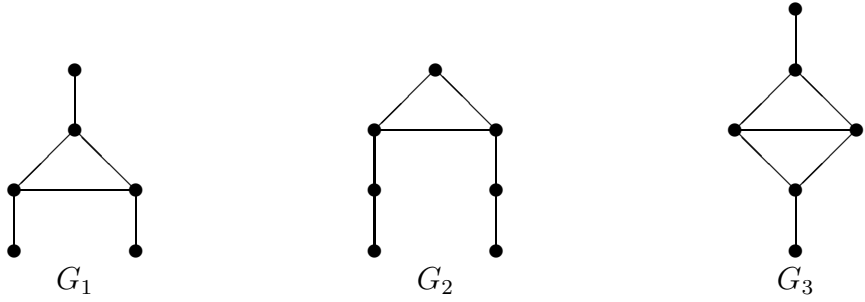

The three graphs $G_1$, $G_2$, and $G_3$ in Figure \ref{fig2} are forbidden subgraphs for the graphs in ${\cal G}_4$.
In fact, each of these graphs contains a $4$-path but has order less than $8$, 
which implies that adding edges yields graphs with $4$-matching number $1$.
Conversely, their $4$-vertex cover number is $2$, and adding edges can only increase this value.

\noindent First, we assume that 
$B$ contains two adjacent vertices $x$ and $y$ 
with exactly $p$ common neighbors $z_1,\ldots,z_p$, where $p\geq 2$.
Let $Z=\{ z_1,\ldots,z_p\}$ and $U=\{ x,y\}\cup Z$. 
If $x$ has a neighbor $x'$ in $B$ outside of $U$, then, since $B$ has no cutvertex, 
a path in $B-x$ between $x'$ and $U\setminus \{ x\}$ together with a suitable path within $B[U]$ yields 
a cycle of order at least $4$ whose order is not a multiple of $4$, contradicting Lemma \ref{lemmacycle}(i).
Hence, $x$, and, by symmetry, $y$ do not have neighbors in $B$ outside of $U$. 
A similar argument also implies that $z_1,\ldots,z_p$ do not have neighbors in $B$ outside of $U$,
which implies that $V(B)=U$.

If $Z$ is not independent, and $p\geq 3$, then $B$ contains a cycle of order $5$, contradicting Lemma \ref{lemmacycle}(i).
Hence, if $Z$ is not independent, then $p=2$, which implies that $B$ is $K_4$.
Since $G$ does not contain $G_1$ as a subgraph, 
we obtain that $B$ is an endblock, that is, 
$B$ is as in (ii) in the definition of ${\cal H}_4$.
Hence, we may assume that $Z$ is independent.
If some vertex in $Z$ is a cutvertex of $G$, 
then, since $G$ does not contain $G_1$ or $G_3$ as a subgraph,
we obtain that $p=2$, and that $B$ is an endblock, that is,
$B$ is as in (iii) in the definition of ${\cal H}_4$.
Hence, we may assume that no vertex in $Z$ is a cutvertex of $G$,
which implies that at most two cutvertices of $G$ belong to $B$,
and that every cutvertex of $G$ that belongs to $B$ has degree $p+1$ in $B$.
Furthermore, if $B$ contains two cutvertices of $G$,
then, since $G$ does not contain $G_2$ as a subgraph,
there is one cutvertex $x$ of $G$ in $B$ such that every vertex in $N_G(x)\setminus V(B)$ has degree $1$ in $G$, that is,
$B$ is as in (iv) in the definition of ${\cal H}_4$.

Next, we assume that $B$ contains a triangle with vertices $x$, $y$, and $z$,
but that no two adjacent vertices in $B$ have more than one common neighbor.
Arguing as above, we obtain $V(B)=\{ x,y,z\}$, and,
since $G$ does not contain $G_1$ or $G_2$ as a subgraph,
it follows that 
$B$ is as in (iv) in the definition of ${\cal H}_4$.
Hence, we may assume that $B$ contains no triangle.

Suppose that $B$ contains a path $P:u_0\ldots u_\ell$ 
such that $u_0$ and $u_\ell$ have degree at least $3$ in $G$,
and $u_1,\ldots,u_{\ell-1}$ have degree $2$ in $G$.
Since $B-u_1$ is connected, 
the path $P$ is contained in a cycle $C$ 
such that $u_0$ and $u_\ell$ both have neighbors outside of $V(C)$.
By Lemma \ref{lemmacycle}(iv), the length $\ell$ of $P$ is even,
in particular, no two vertices of $B$ of degree at least $3$ in $G$ are adjacent.
Let $H$ be the multigraph that arises by replacing every path or cycle
$u_0u_1u_2\ldots u_{2p-2}u_{2p-1}u_{2p}$ of length $2p$ 
such that $u_0$ and $u_{2p}$ have degree at least $3$ in $G$,
and $u_1,\ldots,u_{2p-1}$ have degree $2$ in $G$,
by the path or cycle $u_0u_2\ldots u_{2p-2}u_{2p}$ of length $p$. 
Clearly, $B=Sub_2(H)$, and every cutvertex of $G$ that belongs to $B$ is a vertex of $H$, that is, 
$B$ is as in (i) in the definition of ${\cal H}_4$.
Altogether, it follows that $G\in {\cal H}_4$, which implies ${\cal G}_4\subseteq {\cal H}_4$.

Again, it follows easily from its definition that ${\cal H}_4$ is a hereditary class of graphs.
Hence, in order to show the reverse inclusion ${\cal H}_4\subseteq {\cal G}_4$,
it suffices to show that $\nu_4(G)=\tau_4(G)$ for every connected graph $G$ in ${\cal H}_4$,
which we do by induction on the sum of the order and the size of $G$.
As in the proof of Theorem \ref{theoremg3}, 
we may assume that $G$ is neither a tree nor a cycle.
If $G$ contains a block $B$ as in (ii) or (iii) in the definition of ${\cal H}_4$,
then it is easy to see that $\nu_4(G)=\nu_4(G-V(B))+1$ and $\tau_4(G)=\tau_4(G-V(B))+1$.
If $G$ contains a block $B$ as in (iv) in the definition of ${\cal H}_4$,
then we consider a graph $G'$ obtained from $G$ by removing an edge of $B$ 
that is incident with every cutvertex in $B$.
This graph $G'$ is in ${\cal H}_4$, has less edges than $G$,
and satisfies $\nu_4(G)=\nu_4(G')$ and $\tau_4(G)=\tau_4(G')$.
In all these cases, we obtain $\nu_4(G)=\tau_4(G)$ by induction.
Hence, we may assume that $G$ contains no such block.

Let $B$ be a non-trivial block of $G$.  Let $X$ be the set of
cutvertices of $G$ that belong to $B$.  For $x\in X$, let $G_x$ be the
component of $G-(V(B)\setminus \{ x\})$ that contains $x$.  We may
assume that $B$ is chosen in such a way that there is a vertex $x^*$
in $X$ such that $G_x$ is a tree for every vertex $x$ in $X\setminus
\{ x^*\}$.  If some tree $G_x$ with $x$ in $X\setminus \{ x^*\}$
contains a $4$-path, then Lemma \ref{lemmatree} implies the existence
of an induced subgraph $G'$ of $G$ with $\nu_4(G)=\nu_4(G')+1$ and
$\tau_4(G)=\tau_4(G')+1$, and $\nu_4(G)=\tau_4(G)$ follows by
induction.  Hence, for every vertex $x$ in $X\setminus \{ x^*\}$, the
tree $G_x$ is a star.  Let $X'$ be the set of vertices $x$ in
$X\setminus \{ x^*\}$, for which $G_x$ is not a star with center
vertex $x$, that is, $G_x$ contains a $3$-path $P_x$ starting in $x$.
Let $B'$ be the union of $B$ and the paths $P_x$ for $x$ in $X'$.  If
$B=Sub_2(H)$, where $H$ is as in (i) in the definition of ${\cal
  H}_4$, then $B'=Sub_2(H')$ for the multigraph $H'$ that arises from
$H$ by attaching a vertex of degree $1$ to every vertex in $X'$.
Clearly, $H'$ is bipartite, connected, and contains a cycle.

First, suppose that $x^*$ belongs to some minimum vertex cover in $H'$.
By the K\H{o}nig-Egerv\'ary Theorem, this implies that
every maximum matching in $H'$ contains an edge incident with $x^*$.
Let $M$ be a maximum matching in $H'$.
Similarly as in the proofs of Lemma \ref{lemmasub} and Lemma \ref{lemmabipsub},
we obtain the existence of an injective function $f:M\to E(H')\setminus M$
such that the edges $e$ and $f(e)$ are adjacent for every edge $e$ in $M$. 
Adding the $\nu(H')$ disjoint $4$-paths in $B'$ corresponding to $M$, 
each formed using a subdivided edge $e$ in $M$ and the interior of the subdivided edge $f(e)$,
to a maximum $4$-matching in $G_{x^*}-x^*$ implies $\nu_4(G)\geq \nu_4(G_{x^*}-x^*)+\nu(H')$.
Adding to a minimum $4$-vertex cover in $G_{x^*}-x^*$ 
a minimum vertex cover in $H'$ that contains $x^*$ but none of the vertices of degree $1$ in $V(H')\setminus V(H)$,
yields a $4$-vertex cover in $G$, which implies $\tau_4(G)\leq \tau_4(G_{x^*}-x^*)+\tau(H')$.
Now, by induction and the K\H{o}nig-Egerv\'ary Theorem for $H'$, we obtain
$\nu_4(G)\geq \nu_4(G_{x^*}-x^*)+\nu(H')=\tau_4(G_{x^*}-x^*)+\tau(H')\geq \tau_4(G)\geq \nu_4(G)$,
that is, $\nu_4(G)=\tau_4(G)$.

Now, we may assume that $x^*$ belongs to no minimum vertex cover in $H'$,
which implies that every minimum vertex cover in $H'$ contains all neighbors of $x^*$ in $H'$.
Furthermore, by the K\H{o}nig-Egerv\'ary Theorem,
this implies that some maximum matching $M$ in $H'$ 
contains no edge incident with $x^*$.
Similarly as in the proof of Lemma \ref{lemmabipsub},
we obtain the existence of an injective function
$f:\{ x^*\}\cup M\to E(H')\setminus M$
such that $x^*$ and $f(x^*)$ are incident, 
and $e$ and $f(e)$ are adjacent for every $e\in M$.
Let $G'$ arise from $G_{x^*}$ by attaching a vertex of degree $1$ to $x^*$,
corresponding to the internal vertex of the subdivided version of $f(x^*)$.
Arguing similarly as above,
we obtain $\nu_4(G)\geq \nu_4(G')+\nu(H')$ and $\tau_4(G)\leq \tau_4(G')+\tau(H')$,
and $\nu_4(G)=\tau_4(G)$ follows by induction and the K\H{o}nig-Egerv\'ary Theorem for $H'$,
which completes the proof. $\Box$

\section{Graphs without short cycles in ${\cal G}_k$ for odd $k$}

For general $k$, an explicit characterization of ${\cal G}_k$, similar
to the ones that we obtained for ${\cal G}_3$ and ${\cal G}_4$ in the
previous section, might not be possible. For instance, every graph of
order less than $k$ without a cutvertex is a block of some graph
in ${\cal G}_k$, and already in the characterization of ${\cal G}_4$,
we encountered sporadic blocks that required special attention.
Nevertheless, if we consider an odd $k$ as well as the graphs in
${\cal G}_k$ that do not contain short cycles, then the sporadic
blocks should disappear.

Let $k$ be a positive odd integer.  Let ${\cal G}'_k$ be the set of
all graphs in ${\cal G}_k$ that contain no cycle of order less than
$k$.  Note that ${\cal G}'_3$ actually coincides with ${\cal G}_3$.
Let ${\cal H}'_k$ be the set of all graphs $G$ such that every
non-trivial block $B$ of $G$ satisfies the following condition.
\begin{enumerate}[(i)]
\item $B=Sub_k(H)$ for some multigraph $H$, and every component $K$ of
  $G-V(H)$ that contains a vertex from $V(B)\setminus V(H)$ is a tree
  without a $k$-path.
\end{enumerate}
As before our goal is to show that ${\cal G}'_k$ and ${\cal H}'_k$
coincide.  The following lemma deals with some rather simple graphs in
${\cal G}_k'$ for which it is surprisingly difficult to show that they
belong to ${\cal H}'_k$.

\begin{lemma}\label{lemmaspecial}
Let $k$ be a positive odd integer, and let $p$ be a positive integer.
If the graph $G$ in ${\cal G}_k$ arises from the cycle
$C_{pk}:u_1u_2\ldots u_{pk}u_1$ of order $pk$ by attaching, for every
$i$ in $[pk]$, a path $P_i$ of order $p_i$ to the vertex $u_i$, where
$0\leq p_i<(k-1)/2$, then $G\in {\cal H}'_k$.
\end{lemma}
{\it Proof:} It suffices to show that $\nu_k(G)=p$.
Indeed, if $\nu_k(G)=p$, then $\nu_k(G)=\tau_k(G)=\tau_k(C_{pk})=p$, 
and Lemma \ref{lemmacycle}(ii) implies the existence of a minimum $k$-vertex cover $X$ of $G$ 
with $X=\{ u_{i+(j-1)k}:j\in [p]\}$ for some $i\in [k]$.
It follows that the unique cycle $C_{pk}$ in $G$, which is the only non-trivial block of $G$,
is the $k$-subdivision of the cycle $u_iu_{i+k}u_{i+2k}\ldots u_{i+(p-1)k}u_i$ of order $p$ with vertex set $X$,
and that every component of $G-X$ is a tree without a $k$-path,
that is, $G\in {\cal H}'_k$.
Hence, for a contradiction, we assume that $\nu_k(G)>p$.

Since removing an endvertex from $G$ can reduce the $k$-matching number by at most $1$,
we may assume, by considering a suitable induced subgraph of $G$, 
that $\nu_k(G)=p+1$, and that $\nu_k(G-x)=p$ for every endvertex $x$ of $G$.
For $i$ in $[pk]$, let $P_i$ be the path $u_i^1\ldots u_i^{p_i}$, where, for $p_i\geq 1$, the vertex $u_i^1$ is a neighbor of $u_i$.
Note that the order of $G$ is $pk+p_1+\cdots+p_{pk}$,
and that the endvertices of $G$ are the vertices $u_i^{p_i}$ for those $i$ in $[pk]$ with $p_i\geq 1$. 
Let ${\cal P}$ be a maximum $k$-matching in $G$.
A path $P$ in ${\cal P}$ that is not completely contained in $C_{pk}$ is called {\it special}.
By the choice of $G$, for every special path $P$ in ${\cal P}$, 
there are two distinct indices $i$ and $j$ in $[pk]$ with $\max\{ p_i,p_j\}\geq 1$
such that $P$ is the path 
\begin{eqnarray}\label{e1}
\underbrace{u_i^{p_i}\ldots u_i^1}_{P_i}\underbrace{u_iu_{i+1}\ldots u_{j-1}u_j}_{\subseteq C_{pk}}
\underbrace{u_j^1\ldots u_j^{p_j}}_{P_j},
\end{eqnarray}
where we identify indices modulo $pk$ for the subpath $u_iu_{i+1}\ldots u_{j-1}u_j$ of $P$ that is contained in $C_{pk}$.
If $p_i\geq 1$, then $P$ is said to have the {\it left leg} $P_i$,
If $p_j\geq 1$, then $P$ is said to have the {\it right leg} $P_j$.
Since $G$ contains at most $\nu_k(C_{pk})=p$ disjoint non-special paths, 
and every special path contains at most $2\max\{ p_1,\ldots,p_{pk}\}<k-1$ vertices that do not belong to $C_{pk}$,
the set ${\cal P}$ contains at least two special paths.
By the choice of $G$, for every $i$ in $[pk]$ with $p_i\geq 1$,
the path $P_i$ is either the left leg or the right leg of some path in ${\cal P}$.

Let $i$ in $[pk]$ be such that $P_i$ is the left leg of some path $P$ in ${\cal P}$ as in (\ref{e1}).
By the choice of $G$, the graph $G_i=G-u_i^{p_i}$ satisfies $\nu_k(G_i)=p$.
Similarly as above, 
this implies the existence of a minimum $k$-vertex cover $X_i$ in $G_i$ 
with $X_i=\{ u_{r+(s-1)k}:s\in [p]\}$ for some $r\in [k]$.
We will show that $r=i-p_i$, which implies that $X_i$ is uniquely determined.
Since $X_i$ has order $p$, and intersects all $p$ paths in ${\cal P}\setminus \{ P\}$,
it contains no vertex of $P$, and, hence, no vertex from $u_iu_{i+1}\ldots u_{j-1}u_j$.
Since $G_i-X_i$ contains no $k$-path, this implies that 
$r\in \{ i-p_i,i-p_i+1,\ldots,i-1\}$.
Now, if $r$ is not $i-p_i$, then $r\in \{ i-p_i+1,\ldots,i-1\}$,
the set $X_i$ contains no vertex from $u_iu_{i+1}\ldots u_{i+k-p_i+1}$,
and 
$u_i^{p_i-1}\ldots u_i^1u_iu_{i+1}\ldots u_{i+k-p_i+1}$
is a $k$-path in $G_i-X_i$,
which is a contradiction.
Hence, $r=i-p_i$ as claimed.
Symmetrically, if $P_i$ is the right leg of some path in ${\cal P}$,
then $G_i$ has a unique minimum $k$-vertex cover $X_i$ with $X_i=\{ u_{i+p_i+(s-1)k}:s\in [p]\}$.

We consider some cases.

\medskip

\noindent {\bf Case 1} {\it No path in ${\cal P}$ has a right leg.}

\medskip

\begin{figure}[H]
\begin{center}
{\footnotesize
\unitlength .7mm 
\linethickness{0.4pt}
\ifx\plotpoint\undefined\newsavebox{\plotpoint}\fi 
\begin{picture}(240,18)(0,0)
\put(125,5){\circle*{2}}
\put(125,9){\circle*{2}}
\put(125,13){\circle*{2}}
\thicklines
\put(125,13){\line(0,-1){8}}
\put(125,5){\line(1,0){42}}
\put(167,5){\circle*{2}}
\put(125,0){\makebox(0,0)[cc]{$u_s$}}
\put(171,5){\circle*{2}}
\put(221,5){\circle*{2}}
\put(171,5){\line(1,0){50}}
\put(179,5){\circle*{2}}
\put(225,5){\circle*{2}}
\put(229,5){\circle*{2}}
\put(233,5){\circle*{2}}
\put(225,9){\circle*{2}}
\put(225,17){\circle*{2}}
\put(225,5){\line(0,1){0}}
\put(225,5){\line(1,0){10}}
\put(234.9,4.9){\line(1,0){.8333}}
\put(236.566,4.9){\line(1,0){.8333}}
\put(238.233,4.9){\line(1,0){.8333}}
\put(121,5){\circle*{2}}
\put(71,5){\circle*{2}}
\put(71,5){\line(1,0){50}}
\put(67,5){\circle*{2}}
\put(29,5){\circle*{2}}
\put(29,9){\circle*{2}}
\put(29,17){\circle*{2}}
\put(29,17){\line(0,-1){12}}
\put(29,5){\line(1,0){38}}
\put(29,1){\makebox(0,0)[cc]{$u_r$}}
\put(34,10){\makebox(0,0)[cc]{$u_r^1$}}
\put(37,18){\makebox(0,0)[cc]{$u_r^{p_r-1}$}}
\put(25,5){\circle*{2}}
\put(17,5){\circle*{2}}
\put(13,5){\circle*{2}}
\put(25,5){\line(-1,0){15}}
\put(9,5){\circle*{2}}
\put(9,5){\line(-1,0){3}}
\put(5.9,4.9){\line(-1,0){.8333}}
\put(4.233,4.9){\line(-1,0){.8333}}
\put(2.566,4.9){\line(-1,0){.8333}}
\thinlines
\put(25,5){\line(1,0){4}}
\put(67,5){\line(1,0){4}}
\put(121,5){\line(1,0){4}}
\put(130,10){\makebox(0,0)[cc]{$u_2^1$}}
\put(131,16){\makebox(0,0)[cc]{$u_s^{p_s}$}}
\put(167,5){\line(1,0){4}}
\put(173,7){\vector(4,-3){.1}}\multiput(163,15)(.060240964,-.048192771){166}{\line(1,0){.060240964}}
\put(204,7){\vector(4,-3){.1}}\multiput(194,15)(.060240964,-.048192771){166}{\line(1,0){.060240964}}
\put(155,17){\makebox(0,0)[cc]{$p_s$ vertices}}
\put(193,17){\makebox(0,0)[cc]{$(k-p_s-1)$ vertices}}
\put(225,0){\makebox(0,0)[cc]{$u_t$}}
\put(229,10){\makebox(0,0)[cc]{$u_t^1$}}
\put(230,18){\makebox(0,0)[cc]{$u_t^{p_t}$}}
\put(11,3){\framebox(4,4)[cc]{}}
\put(69,3){\framebox(4,4)[cc]{}}
\put(123,3){\framebox(4,4)[cc]{}}
\put(172.5,5){\oval(7,4)[]}
\put(177,3){\framebox(4,4)[cc]{}}
\put(202.5,5){\oval(41,4)[]}
\put(231,3){\framebox(4,4)[cc]{}}
\put(221,5){\line(1,0){4}}
\thicklines
\put(225,17){\line(0,-1){12}}
\end{picture}
}
\end{center}
\caption{The situation in Case 1, where vertices in $X_r$ are indicated by the square boxes, and the paths in ${\cal P}$ are shown in bold.}
\label{fig:Case1}
\end{figure}
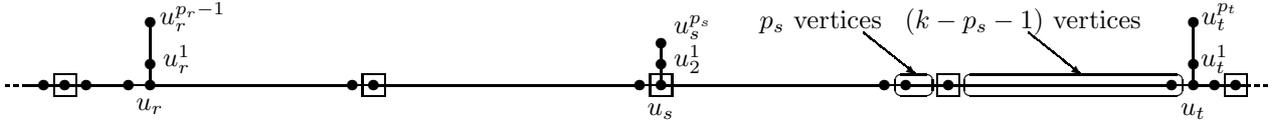

\noindent In this case, every special path in ${\cal P}$ contains at most $\max\{ p_1,\ldots,p_{pk}\}<(k-1)/2$ 
vertices that do not belong to $C_{pk}$,
which implies that ${\cal P}$ contains at least three special paths.
By symmetry, we may assume that the indices $r$, $s$, and $t$ in $[pk]$ are chosen in such a way that 
\begin{itemize}
\item $r<s<t$,
\item $p_s\le p_t$,
\item $P_r$, $P_s$, and $P_t$ are left legs of three special paths in ${\cal P}$, and 
\item no other special path in ${\cal P}$ intersects the subpath $u_r\ldots u_s\ldots u_t$ of $C_{pk}$.
\end{itemize}
By the choice of $G$, in this case it follows that every vertex of $C_{pk}$ belongs to some path in ${\cal P}$.
Therefore, the final condition in the choice of $r$, $s$, and $t$ implies that 
$$s\equiv (r+k-p_r)\mod k\,\,\,\,\,\mbox{ and }\,\,\,\,\,t\equiv
(s+k-p_s)\mod k.$$ Since $X_r$ contains the vertex $u_{r-p_r}$, this
implies that $u_s\in X_r$, and that $X_r$ contains no vertex from
$u_{t-k+p_s+1}u_{t-k+p_s+2}\ldots u_t$. 
See Figure~\ref{fig:Case1} for an illustration.

\noindent Nevertheless, since $p_s\leq p_t$, the graph $G_r-X_r$ contains the
path $u_{t-k+p_s+1}u_{t-k+p_s+2}\ldots u_tu_t^1\ldots u_t^{p_t}$ of
order $k-1-p_s+1+p_t\geq k$, which is a contradiction.

\medskip

\noindent {\bf Case 2} {\it Some special path in ${\cal P}$ has a right leg, and some special path in ${\cal P}$ has a left leg.}

\medskip

\begin{figure}[H]
\begin{center}
{\footnotesize
\unitlength 1mm 
\linethickness{0.4pt}
\ifx\plotpoint\undefined\newsavebox{\plotpoint}\fi 
\begin{picture}(155,18)(0,0)
\put(86,5){\circle*{2}}
\put(136,5){\circle*{2}}
\thicklines
\put(86,5){\line(1,0){50}}
\put(94,5){\circle*{2}}
\put(140,5){\circle*{2}}
\put(144,5){\circle*{2}}
\put(148,5){\circle*{2}}
\put(140,9){\circle*{2}}
\put(140,17){\circle*{2}}
\put(140,5){\line(0,1){0}}
\put(140,5){\line(1,0){10}}
\put(149.93,4.93){\line(1,0){.8333}}
\put(151.596,4.93){\line(1,0){.8333}}
\put(153.263,4.93){\line(1,0){.8333}}
\put(82,5){\circle*{2}}
\put(32,5){\circle*{2}}
\put(32,5){\line(1,0){50}}
\thinlines
\put(82,5){\line(1,0){4}}
\put(88,7){\vector(4,-3){.07}}\multiput(78,15)(.042016807,-.033613445){238}{\line(1,0){.042016807}}
\put(119,7){\vector(4,-3){.07}}\multiput(109,15)(.042016807,-.033613445){238}{\line(1,0){.042016807}}
\put(70,17){\makebox(0,0)[cc]{$(p_s-s'-1)$ vertices}}
\put(108,17){\makebox(0,0)[cc]{$(k-p_s+s')$ vertices}}
\put(140,0){\makebox(0,0)[cc]{$u_t$}}
\put(144,10){\makebox(0,0)[cc]{$u_t^1$}}
\put(145,18){\makebox(0,0)[cc]{$u_t^{p_t}$}}
\put(38,3){\framebox(4,4)[cc]{}}
\put(87.5,5){\oval(7,4)[]}
\put(92,3){\framebox(4,4)[cc]{}}
\put(117.5,5){\oval(41,4)[]}
\put(146,3){\framebox(4,4)[cc]{}}
\put(136,5){\line(1,0){4}}
\put(40,5){\circle*{2}}
\put(28,5){\circle*{2}}
\put(20,5){\circle*{2}}
\put(16,5){\circle*{2}}
\put(16,9){\circle*{2}}
\put(16,17){\circle*{2}}
\put(16,5){\line(1,0){16}}
\put(28,1){\makebox(0,0)[cc]{$u_{s+s'}$}}
\put(16,1){\makebox(0,0)[cc]{$u_s$}}
\put(20,10){\makebox(0,0)[cc]{$u_s^1$}}
\put(23,18){\makebox(0,0)[cc]{$u_s^{p_s}-1$}}
\thicklines
\put(16,17){\line(0,-1){12}}
\put(16,5){\line(-1,0){12}}
\put(3.93,4.93){\line(-1,0){.8333}}
\put(2.263,4.93){\line(-1,0){.8333}}
\put(.596,4.93){\line(-1,0){.8333}}
\put(140,17){\line(0,-1){12}}
\end{picture}
}
\end{center}
\caption{Illustration of the proof that $s'=p_s-1$.}
\label{fig:Case2a}
\end{figure}
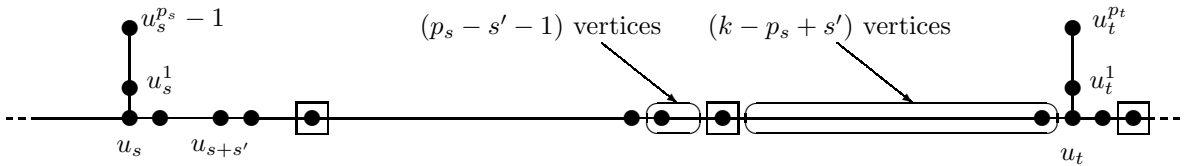

\noindent By symmetry, we may assume that the indices $s$ and $t$ in $[pk]$ are such that 
\begin{itemize}
\item $s<t$,
\item $p_s\leq p_t$,
\item $P_s$ is the right leg of a special path in ${\cal P}$, and $P_t$ is the left leg of a special path in ${\cal P}$, and 
\item no other special path in ${\cal P}$ intersects the subpath $u_s\ldots u_t$ of $C_{pk}$.
\end{itemize}
We may assume that the non-special paths in ${\cal P}$ that intersect
$u_s\ldots u_t$ are chosen in such a way that their removal from
$u_s\ldots u_t$ leaves a path of the form $u_s\ldots u_{s+s'}$ for
some $s'\geq 0$.  Since $\nu_k(G_s)=p$, we have $s'\leq p_s-1$.  If
$s'\leq p_s-2$, then $X_s$ contains no vertex from
$u_{t-k+p_s-s'}u_{t-k+p_s-s'+1}\ldots u_t$, and $G_s-X_s$ contains the
path $u_{t-k+p_s-s'}u_{t-k+p_s-s'+1}\ldots u_t u_t^1\ldots u_t^{p_t}$
of order $k-p_s+s'+1+p_t>k$, which is a contradiction. 
See Figure~\ref{fig:Case2a} for an illustration.
Hence, we obtain
$s'=p_s-1$, which implies that $u_t\in X_s$.

\begin{figure}[H]
\begin{center}
{\footnotesize
\unitlength .8mm 
\linethickness{0.4pt}
\ifx\plotpoint\undefined\newsavebox{\plotpoint}\fi 
\begin{picture}(186,18)(0,0)
\put(120,5){\circle*{2}}
\put(170,5){\circle*{2}}
\thicklines
\put(120,5){\line(1,0){50}}
\put(124,5){\circle*{2}}
\put(174,5){\circle*{2}}
\put(178,5){\circle*{2}}
\put(174,9){\circle*{2}}
\put(174,17){\circle*{2}}
\put(174,5){\line(0,1){0}}
\put(116,5){\circle*{2}}
\put(66,5){\circle*{2}}
\put(66,5){\line(1,0){50}}
\thinlines
\put(116,5){\line(1,0){4}}
\put(174,0){\makebox(0,0)[cc]{$u_t$}}
\put(178,10){\makebox(0,0)[cc]{$u_t^1$}}
\put(179,18){\makebox(0,0)[cc]{$u_t^{p_t}$}}
\put(68,3){\framebox(4,4)[cc]{}}
\put(122,3){\framebox(4,4)[cc]{}}
\put(176,3){\framebox(4,4)[cc]{}}
\put(170,5){\line(1,0){4}}
\put(70,5){\circle*{2}}
\put(62,5){\circle*{2}}
\put(54,5){\circle*{2}}
\put(50,5){\circle*{2}}
\put(50,9){\circle*{2}}
\put(50,17){\circle*{2}}
\put(50,5){\line(1,0){16}}
\put(50,1){\makebox(0,0)[cc]{$u_s$}}
\put(54,10){\makebox(0,0)[cc]{$u_s^1$}}
\put(55,18){\makebox(0,0)[cc]{$u_s^{p_s}$}}
\thicklines
\put(50,17){\line(0,-1){12}}
\put(70,0){\makebox(0,0)[cc]{$u_{s+p_s+1}$}}
\put(50,5){\line(-1,0){30}}
\put(20,5){\line(0,1){8}}
\put(20,5){\circle*{2}}
\put(20,9){\circle*{2}}
\put(20,13){\circle*{2}}
\put(20,1){\makebox(0,0)[cc]{$u_r$}}
\put(24,10){\makebox(0,0)[cc]{$u_r^1$}}
\put(27,16){\makebox(0,0)[cc]{$u_r^{p_r-1}$}}
\put(16,5){\circle*{2}}
\put(8,5){\circle*{2}}
\thinlines
\put(20,5){\line(-1,0){12}}
\put(3.912,4.912){\line(-1,0){.8333}}
\put(2.245,4.912){\line(-1,0){.8333}}
\put(.579,4.912){\line(-1,0){.8333}}
\put(12,5){\circle*{2}}
\thicklines
\put(174,5){\line(1,0){7}}
\put(180.912,4.912){\line(1,0){.8333}}
\put(182.579,4.912){\line(1,0){.8333}}
\put(184.245,4.912){\line(1,0){.8333}}
\thinlines
\put(5,5){\line(1,0){3}}
\put(6,3){\framebox(4,4)[cc]{}}
\thicklines
\put(174,17){\line(0,-1){12}}
\end{picture}
}
\end{center}
\caption{Illustration of the proof that $P'$ has no left leg.}
\label{fig:Case2b}
\end{figure}
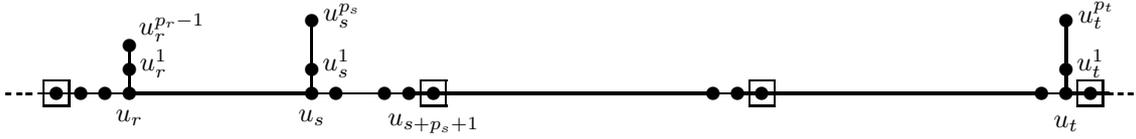

\noindent If the path $P'$ in ${\cal P}$ whose right leg is $P_s$ also has a
left leg, say $P_r$ for some $r<s$, then $X_r$ contains $x_{r-p_r}$,
and, hence, also $u_{s+p_s+1}$ as well as $u_{t+1}$ but no vertex from
$u_{t-k+2}u_{t-k+3}\ldots u_t$. 
See Figure~\ref{fig:Case2b} for an illustration.

\noindent Now, $G_r-X_r$ contains the path
$u_{t-k+2}u_{t-k+3}\ldots u_t u_t^1\ldots u_t^{p_t}$ of order
$k-2+1+p_t\geq k$, which is a contradiction.  Hence, $P'$ has no left
leg, and equals $u_{s-k+p_s+1}u_{s-k+p_s+2}\ldots u_s u_s^1\ldots
u_s^{p_s}$.

Let $r<s$ be maximum such that some special path $P''$ in ${\cal P}$
contains $u_r$.  By the choice of $G$, and, since $P'$ has no left
leg, we obtain that $r\equiv (s-k+p_s)\mod k$.

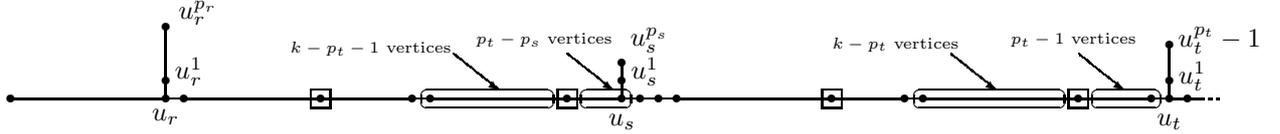
\begin{figure}[H]
\begin{center}
{\footnotesize
\unitlength .6mm 
\linethickness{0.4pt}
\ifx\plotpoint\undefined\newsavebox{\plotpoint}\fi 
\begin{picture}(266,22)(0,0)
\put(134,5){\circle*{2}}
\put(134,9){\circle*{2}}
\put(134,13){\circle*{2}}
\thicklines
\put(134,13){\line(0,-1){8}}
\put(134,0){\makebox(0,0)[cc]{$u_s$}}
\put(200,5){\circle*{2}}
\put(146,5){\circle*{2}}
\put(250,5){\circle*{2}}
\put(196,5){\circle*{2}}
\put(234,5){\circle*{2}}
\put(180,5){\circle*{2}}
\put(200,5){\line(1,0){50}}
\put(146,5){\line(1,0){50}}
\put(254,5){\circle*{2}}
\put(258,5){\circle*{2}}
\put(254,9){\circle*{2}}
\put(254,17){\circle*{2}}
\put(254,5){\line(0,1){0}}
\put(139,11){\makebox(0,0)[cc]{$u_s^1$}}
\put(140,18){\makebox(0,0)[cc]{$u_s^{p_s}$}}
\thinlines
\put(151,5){\line(1,0){4}}
\put(216,7){\vector(4,-3){.117}}\multiput(206,15)(.06993007,-.055944056){143}{\line(1,0){.06993007}}
\put(244,7){\vector(4,-3){.117}}\multiput(234,15)(.06993007,-.055944056){143}{\line(1,0){.06993007}}
\put(132,7){\vector(4,-3){.117}}\multiput(122,15)(.06993007,-.055944056){143}{\line(1,0){.06993007}}
\put(107,7){\vector(4,-3){.117}}\multiput(97,15)(.06993007,-.055944056){143}{\line(1,0){.06993007}}
\put(194,17){\makebox(0,0)[cc]{\tiny $k-p_t$ vertices}}
\put(233,18){\makebox(0,0)[cc]{\tiny $p_t-1$ vertices}}
\put(117,18){\makebox(0,0)[cc]{\tiny $p_t-p_s$ vertices}}
\put(79,16){\makebox(0,0)[cc]{\tiny $k-p_t-1$ vertices}}
\put(254,0){\makebox(0,0)[cc]{$u_t$}}
\put(259,10){\makebox(0,0)[cc]{$u_t^1$}}
\put(265,18){\makebox(0,0)[cc]{$u_t^{p_t}-1$}}
\put(250,5){\line(1,0){4}}
\put(196,5){\line(1,0){4}}
\thicklines
\put(254,17){\line(0,-1){12}}
\put(254,5){\line(1,0){7}}
\put(260.883,4.883){\line(1,0){.8333}}
\put(262.55,4.883){\line(1,0){.8333}}
\put(264.216,4.883){\line(1,0){.8333}}
\put(138,5){\circle*{2}}
\put(142,5){\circle*{2}}
\thinlines
\put(146,5){\line(-1,0){12}}
\thicklines
\put(134,5){\line(-1,0){42}}
\put(122,5){\circle*{2}}
\put(92,5){\circle*{2}}
\put(88,5){\circle*{2}}
\put(38,5){\circle*{2}}
\put(68,5){\circle*{2}}
\put(38,5){\line(1,0){49}}
\thinlines
\put(66,3){\framebox(4,4)[cc]{}}
\put(91,5){\line(-1,0){3}}
\put(34,5){\circle*{2}}
\put(34,9){\circle*{2}}
\put(34,21){\circle*{2}}
\thicklines
\put(34,21){\line(0,-1){16}}
\put(34,5){\line(-1,0){34}}
\put(0,5){\circle*{2}}
\put(34,1){\makebox(0,0)[cc]{$u_r$}}
\put(39,11){\makebox(0,0)[cc]{$u_r^1$}}
\put(41,24){\makebox(0,0)[cc]{$u_r^{p_r}$}}
\thinlines
\put(34,5){\line(1,0){4}}
\put(104.5,5){\oval(29,4)[]}
\put(130.5,5){\oval(11,4)[]}
\put(120,3){\framebox(4,4)[cc]{}}
\put(214.5,5){\oval(33,4)[]}
\put(232,3){\framebox(4,4)[cc]{}}
\put(178,3){\framebox(4,4)[cc]{}}
\put(244.5,5){\oval(15,4)[]}
\end{picture}
}
\end{center}
\caption{Illustration of the proof that $p_r\not\geq p_t$.}
\label{fig:Case2c}
\end{figure}

First, suppose that $p_r=0$, that is, $P''$ has no right leg. Since
$P''$ is special, it has a left leg, say $P_q$ for some $q<r$. 
Here things work as previously; 
$X_q$ contains $u_{q-p_q}$, and, hence,
also $u_{r+1}$, $u_{s-k+p_s+1}$, $u_{s+p_s+1}$, and $u_{t+1}$ but no
vertex from $u_{t-k+2}u_{t-k+3}\ldots u_t$.  Now, $G_q-X_q$ contains
the path $u_{t-k+2}u_{t-k+3}\ldots u_t u_t^1\ldots u_t^{p_t}$ of order
$k-2+1+p_t\geq k$, which is a contradiction.
Hence, $p_r\geq 1$, that is, the path $P''$ has $P_r$ as its
right leg.  If $p_r\geq p_t$, then $X_t$ contains $u_{t-p_t}$, and,
hence, also $u_{s-p_t+p_s}$ as well as $u_{r+k-p_t}$ but no vertex
from $u_ru_{r+1}\ldots u_{r+k-p_t-1}$. 
See Figure~\ref{fig:Case2c} for an illustration.

\noindent Now, $G_t-X_t$ contains the
path $u_r^{p_r}\ldots u_r^1 u_r u_{r+1}\ldots u_{r+k-p_t-1}$ of order
$p_r+1+k-p_t-1\geq k$, which is a contradiction.  Conversely, if
$p_r<p_t$, then $X_r$ contains $u_{r+p_r}$, and, hence, also
$u_{t-k+p_r}$ but no vertex from $u_{t-k+p_r+1}u_{t-k+p_r+2}\ldots
u_t$.  Now, $G_r-X_r$ contains the path
$u_{t-k+p_r+1}u_{t-k+p_r+2}\ldots u_t u_t^1\ldots u_t^{p_t}$ of order
$k-p_r-1+1+p_t\geq k$, which is a contradiction. This completes the
proof. $\Box$

\medskip

\noindent We proceed to the main result in this section,
which actually contains Theorem \ref{theoremg3} as a special case.
In view of its simplicity, we kept the separate proof of Theorem \ref{theoremg3}.

\begin{theorem}\label{theoremgk}
${\cal G}'_k={\cal H}'_k$ for every positive odd integer $k$.
\end{theorem}
{\it Proof:} As before, in order to show that ${\cal G}'_k\subseteq {\cal H}'_k$,
we show that $G\in {\cal H}'_k$ for every connected graph $G\in {\cal G}'_k$.
By Lemma \ref{lemmacycle}(i), the order of every cycle in $G$ is a multiple of $k$.
We may again assume that $G$ is neither a tree nor a cycle.
Let $B$ be a non-trivial block of $G$.

First, we assume that $B$ is not just a cycle, that is, it contains vertices that are of degree at least $3$ in $B$.
Suppose that $B$ contains a path $P:u_0\ldots u_\ell$ 
such that $u_0$ and $u_\ell$ have degree at least $3$ in $B$,
and $u_1,\ldots,u_{\ell-1}$ have degree $2$ in $B$.
Since $B-u_1$ is connected, 
the path $P$ is contained in a cycle $C$ 
such that $u_0$ and $u_\ell$ both have neighbors outside of $V(C)$, say $u_0^1$ and $u_\ell^1$, respectively.
Let $P_0$ be a shortest path in $B-u_0$ between $u_0^1$ and $V(C)\setminus \{ u_0\}$.
Since the order of every cycle in $G$ is a multiple of $k$, and, since $k$ is odd, 
it follows that $P_0$ has length $-1$ modulo $k$, 
which implies that $B-V(C)$ contains a path $P_0'$ of order $(k-1)/2$ starting in $u_0^1$.
Similarly, $B-V(C)$ contains a path $P'_\ell$ of order $(k-1)/2$ starting in $u_\ell^1$.
If $G'$ is the subgraph of $G$ induced by $V(C)\cup V(P_0')\cup V(P'_\ell)$,
then $\frac{n(C)}{k}\leq \nu_k(G')\leq \left\lfloor\frac{n(C)+n(P_0')+n(P'_\ell)}{k}\right\rfloor=\frac{n(C)}{k}$.
It follows that every minimum $k$-vertex cover $X'$ of $G'$
is also a minimum $k$-vertex cover of $C$, and, hence, as described in Lemma \ref{lemmacycle}(ii).
In view of $P_0'$, $P'_\ell$, and the subpaths of $C$ not covered by $X'$, 
it follows that the vertices $u_0$ and $u_\ell$ must both belong to $X'$.
This implies that the length $\ell$ of $P$ is a multiple of $k$. 
Let $H$ be the multigraph that arises by replacing every path or cycle
$u_0u_1\ldots u_{pk}$ of length $pk$ 
such that $u_0$ and $u_{pk}$ have degree at least $3$ in $B$,
and $u_1,\ldots,u_{pk-1}$ have degree $2$ in $B$,
by the path or cycle $u_0u_k\ldots u_{pk}$ of length $p$. 
Clearly, $B=Sub_k(H)$.

Let $K$ be a component of $G-V(H)$ that contains a vertex from $V(B)\setminus V(H)$.
Let $uv$ be an edge of $H$ such that $K$ intersects the subdivided edge $uv$.
Since $B$ is a block of $G$, the component $K$ intersects $V(B)\setminus V(H)$ exactly in the interior of the subdivided edge $uv$.
Let $P:uw_1\ldots w_{k-1}v$ be the path in $G$ corresponding to the subdivided edge $uv$.
Suppose, for a contradiction, that $K$ contains a $k$-path.
This implies that we may assume, by symmetry, that there is some $i\in [(k-1)/2]$,
and a path $Q:x_1\ldots x_i$ in $K-V(B)$ such that $x_i$ is adjacent to $w_i$.
Let $C$ be a cycle in $B$ containing $P$.
Similarly as above, we obtain the existence of a path $R$ of order $(k-1)/2$ in $B-V(C)$ 
such that $u$ is adjacent to an endvertex of $R$.
If $G'$ is the subgraph of $G$ induced by $V(C)\cup V(Q)\cup V(R)$,
then $\nu_k(G')=\frac{n(C)}{k}$.
Therefore, every minimum $k$-vertex cover of $G'$ is also a minimum $k$-vertex cover of $C$,
and, hence, as described in Lemma \ref{lemmacycle}(ii).
In view of $R$ and the subpaths of $C$ not covered by $X'$, 
it follows that $u$ must belong to $X'$.
But now, $x_1\ldots x_i w_i\ldots w_{k-1}$ is a $k$-path in $G'-X'$,
which is a contradiction.
Altogether, it follows that $K$ contains no $k$-path,
which implies that $K$ is a tree without a $k$-path.
Hence, $B$ is as in (i) in the definition of ${\cal H}'_k$.

Next, we assume that $B$ is a cycle $C:u_1\ldots u_{pk}$.
For every $i$ in $[pk]$, let $p_i$ be the maximum length of a path in $G-(V(B)\setminus \{ u_i\})$ starting in the vertex $u_i$.
First, suppose that $\max\{ p_1,\ldots,p_{pk}\}\geq (k-1)/2$.
By symmetry, we may assume that $p_1\geq (k-1)/2$.
Let $X=\{ u_{1+(j-1)k}:j\in [p]\}$.
Clearly, $B=Sub_k(H)$, where $H$ is the cycle $u_1u_{1+k}\ldots u_{1+(p-1)k}u_1$ with vertex set $X$.

Let $K$ be a component of $G-V(H)$ that contains a vertex from $V(B)\setminus V(H)$.
If $K$ contains a $k$-path,
then, by symmetry, we may assume that there is some index $i$ in $[pk]$ such that $1\leq (i-1)\mod k\leq (k-1)/2$
and, $p_i$ is at least $(i-1)\mod k$.
Now, $G$ contains a subgraph $G'$ that arises from $B$ by attaching a path of order $(k-1)/2$ to $u_1$,
and a path of order $(i-1)\mod k$ to $u_i$.
As before $\nu_k(G')=\frac{n(B)}{k}$, and Lemma \ref{lemmacycle}(ii) implies 
that every minimum $k$-vertex cover $X'$ of $G'$ must contain $u_1$,
and that $G'-X'$ still contains a $k$-path using the path attached to $u_i$, 
which is a contradiction.
Altogether, it follows that $K$ contains no $k$-path,
which implies that $K$ is a tree without a $k$-path.
Hence, $B$ is as in (i) in the definition of ${\cal H}'_k$.

Now, we may assume that $\max\{ p_1,\ldots,p_{pk}\}<(k-1)/2$.
This implies that, for every $i$ in $[pk]$,
the component $G_{u_i}$ of $G-(V(B)\setminus \{ u_i\})$ that contains $u_i$,
is a tree without a $k$-path.
Let $G'$ be the induced subgraph of $G$ that arises from $G$ by removing,
for every $i$ in $[pk]$, 
all of $G_{u_i}$ except for a path of length $p_i$ starting in the vertex $u_i$.
By Lemma \ref{lemmaspecial}, the graph $G'$ belong to ${\cal H}'_k$,
which easily implies that also $G$ belongs to ${\cal H}'_k$.
Altogether, we obtain ${\cal G}'_k\subseteq {\cal H}'_k$.

Again, it follows easily from its definition that ${\cal H}'_k$ is a hereditary class of graphs, and, hence,
in order to show the reverse inclusion ${\cal H}'_k\subseteq {\cal G}'_k$,
it suffices to show that $\nu_k(G)=\tau_k(G)$ for every connected graph $G$ in ${\cal H}'_k$.
This now follows very easily by induction on the order using Lemma \ref{lemmatree} and Lemma \ref{lemmasub},
which completes the proof.
$\Box$

\section{Conclusion}

It is not difficult to extract from our results 
all minimal forbidden induced subgraphs 
for the graph classes ${\cal G}_3$, ${\cal G}_4$, and ${\cal G}'_k$ for odd $k$ at least $5$.
Furthermore, our results imply that the graphs in these classes can be recognized efficiently,
and that there are simple combinatorial polynomial time algorithms 
that determine maximum $k$-matchings and minimum $k$-vertex covers
for these graphs. Apart from extending our characterizations, 
a natural open problem concerns the complexity of recognizing the graphs in ${\cal G}_k$ for general fixed $k$. 
We pose the following optimistic conjecture.

\begin{conjecture}\label{conjecture1}
For every fixed positive integer $k$, it can be decided in polynomial time
whether a given graph belongs to ${\cal G}_k$.
\end{conjecture}
Lemma \ref{lemmacycle}(i) easily implies that every graph in ${\cal G}_k$ has minimum degree at most $k$.
This implies that the graphs in ${\cal G}_k$ are $k$-degenerate, which might be a useful property for their recognition.

For $k\in \{ 3,4\}$, our results imply that $\nu_k(H)=\tau_k(H)$
for every not necessarily induced subgraph $H$ of every graph $G$ in ${\cal G}_k$.
For $k=1$, the same trivially holds, and, also for $k=2$, 
the same holds, since graphs are bipartite if and only if
all their not necessarily induced subgraphs are bipartite.
We believe that these observations generalize, and pose the following conjecture.

\begin{conjecture}\label{conjecture2}
For every positive integer $k$,
the set ${\cal G}_k$ equals the set of all graphs $G$
such that $\nu_k(H)=\tau_k(H)$ for every subgraph $H$ of $G$. 
\end{conjecture}
One proof of the K\H{o}nig-Egerv\'ary Theorem, as well as many polyhedral insights concerning matchings in bipartite graphs,
rely on the total unimodularity of the vertex versus edge incidence matrices of bipartite graphs.
Unfortunately, for integers $k$ at least $3$, the vertex versus $k$-path incidence matrices of the graphs in ${\cal G}_k$ are not totally unimodular.
If $G=Sub_3(H)$ for some graph $H$ with a vertex $u$ of degree at least $3$ for instance,
then considering three suitable $3$-paths containing $u$ as central vertex, and three suitable neighbors of $u$ on these paths, 
implies that the vertex versus $3$-path incidence matrix $A$ of $G$
contains the vertex versus edge incidence matrix of $C_3$ as a submatrix, that is, $A$ is not totally unimodular.


\begin{thebibliography}{}
\bibitem{bcl} R. Boliac, K. Cameron, V.V. Lozin, On computing the dissociation number and the induced matching number of bipartite graphs, Ars Combinatorics 72 (2004) 241-253.
\bibitem{bjkst} B. Bre\v{s}ar, M. Jakovac, J. Katreni\v{c}, G. Semani\v{s}in, A. Taranenko, On the vertex $k$-path cover, Discrete Applied Mathematics 161 (2013) 1943-1949.
\bibitem{bkks} B. Bre\v{s}ar, F. Kardo\v{s}, J. Katreni\v{c}, G. Semani\v{s}in, Minimum $k$-path vertex cover, Discrete Applied Mathematics 159 (2011) 1189-1195.
\bibitem{e} E. Egerv\'ary,  \"Uber kombinatorische Eigenschaften von Matrizen, Matematikai \'es Fizikai Lapok 38 (1931) 16-28.
\bibitem{gj} M.R. Garey, D.S. Johnson, Computers and Intractability: A Guide to the Theory of NP-Completeness, W.H. Freeman and Co. pp. x+338, 1979. 
\bibitem{kks} F. Kardo\v{s}, J. Katreni\v{c}, I. Schiermeyer, On computing the minimum $3$-path vertex cover and dissociation number of graphs, Theoretical Computer Science 412 (2011) 7009-7017.
\bibitem{k2} J. Katreni\v{c}, A faster FPT algorithm for $3$-path vertex cover, Information Processing Letters 116 (2016) 273-278.
\bibitem{k} D. K\H{o}nig, Graphen und Matrices, Matematikai \'es Fizikai Lapok 38 (1931) 116-119.
\bibitem{lzh} X. Li, Z. Zhang, X. Huang, Approximation algorithms for minimum (weight) connected $k$-path vertex cover, Discrete Applied Mathematics 205 (2016) 101-108.
\bibitem{twyc} J. Tu, L. Wu, J. Yuan, L. Cui, On the vertex cover $P_3$ problem parameterized by treewidth, Journal of Combinatorial Optimization 34 (2017) 414-4
\bibitem{tz} J. Tu, W. Zhou, A primal-dual approximation algorithm for the vertex cover $P_3$ problem, Theoretical Computer Science 412 (2011) 7044-7048.
\bibitem{xk} M. Xiao, S. Kou, Exact algorithms for the maximum dissociation set and minimum $3$-path vertex cover problems, Theoretical Computer Science 657 (2017) 86-97.
\bibitem{y} M. Yannakakis, Node-deletion problems on bipartite graphs, SIAM Journal on Computing 10 (1981) 310-327.
\end{thebibliography}
\end{document}